\magnification\magstep1
\input AHTOH-E.STY
\hfuzz 6 pt
\def\KOO{{\rm KOO}(G_1\!)}
\def\Gn{G(x_1,\dots,x_n)}

\UDC{
512.543.72
%
+512.544.33
+512.554.33
}


\title{Solvability of unimodular equations in groups and Lie algebras}

\author{
Anton A. Klyachko$^\sharp$
\quad
Mikhail A. Mikheenko$^\sharp$
\quad
Alexander Yu. Olshanskii$^{\flat\sharp}$
}
\address{
$^\sharp$%
Faculty of Mechanics and Mathematics of Moscow State University
\\
Moscow 119991, Leninskie gory, MSU.
\\
$^\flat$%
Department of Mathematics,
Vanderbilt University, Nashville, TN 37240, U.S.A.
\\
klyachko@mech.math.msu.su
\quad
mamikheenko@mail.ru
\quad
alexander.olshanskiy@vanderbilt.edu
}
\grantsFirstSecond{\RSF 22-11-00075}

\abstract{%
Our results imply, in particular, that a finitely generated solvable
group~$G$ is nilpotent if and only if it contains a solution to any
unimodular equation, i.e., an equation of the form~$\prod g_ix^{n_i}=1$,
where~$g_i\in G$ and $\sum n_i=\pm1$.  A similar fact turns out to be true
for Lie algebras.  We also exhibit an example of a unimodular equation
$w(x)=g$ over a finitely generated group $G$, which has a solution (in
$G$) for any~$g\in G$, but  the solution is not unique for some $g\in G$.
We show that, for nilpotent groups $G$, the set of unimodular mappings
$G^n\to G^n$ (which are defined naturally) forms a group under
the composition.
}


\s 1.
Introduction

An equation $g_1x^{n_1}\dots g_kx^{n_k}=1$ over a group $G$
(where $g_i\in G$ and $n_i\in\Z$) is called \emph{nonsingular} if
$\sum n_i\ne0$; if $\sum n_i=\pm1$, then the equation is called
\emph{unimodular}; the left-hand sides of such equations are called
nonsingular [unimodular] words.  These definitions extend naturally to any
systems of equations with any sets of unknowns (see, e.g., [KMR24]),
however, this paper largely deals with the case of one equation with one
unknown. The only exception is Section 4, where we consider
\emph{unimodular systems} of $n$ equations with $n$ unknowns, i.e., such
systems of equations over groups that the determinant of the matrix 
composed of the exponent sums of the $i$th unknown in the $j$th equation 
is $\pm1$.

\enditem
Generally,
\-
nonsingular equations behave better than singular ones;
e.g., the singular equation $x^{-1}axb=1$,
where $a,b\in G$ are elements of different orders,
has no solutions in any overgroup of $G$;
the existence of a nonsingular equation
with such property is a well-known open question ---
the Kervaire--Laudenbach conjecture%
\fn{Quite resently, \emph{OpenAI} claimed these
questions solved [OAI26a], [OAI26b].}%
;

\-
unimodular equations behave better than arbitrary nonsingular ones;
e.g.,
the (nonsingular but not unimodular)
equation $x^2g=1$ has no solutions
in the infinite cyclic group $\gp g$;
but no unimodular equation with such property
exists over
nilpotent groups (see below);
another example:
in~[K93] (see also [FeR96]), it was shown that
\disp{\sl
any unimodular equation over torsion-free group has a solution in
some overgroup;
}%
while it is unknown whether a similar assertion holds for
arbitrary nonsingular equations$^{*)}$.

\proclaim Shmelkin theorem \rm[Sh67] (a simplified and weakened version).
Any finite unimodular system of equations over a nilpotent group $G$ has a
unique solution in $G$.

\noindent
For solvable groups, the situation is more complicated:
\-
for each $n>1$,
there exists a unimodular equation
with one unknown
over a solvable group of derived length $n$,
which has no solutions in any larger
solvable group of the same derived
length $n$ [M26b];
\-
but,
for each positive integer $n$,
any nonsingular system equations
over any group~$G$
with a subnormal
series
$
G=G_1\triangleright\dots\triangleright G_{n+1}=\1,
$
where all sections $G_i/G_{i+1}$ are torsion-free abelian,
has a solution in an overgroup having a subnormal series
of the same length $n$, whose all sections are torsion-free abelian
[KMR24].

\enditem
In Section 3, we supplement these results from [Sh67], [KMR24], and [M26b]
with the following proposition:
\disp{\sl
Over each finitely generated solvable \(or linear\)
non-nilpotent group,
there exists a unimodular equation having no solutions in this group.
}%
We suggest to call a group $G$ \emph{unimodularly \[nonsingularly\]
closed} if any unimodular [nonsingular] equation over $G$ has a solution
in~$G$.
Thus,
a finitely generated solvable \(or linear\) group is unimodularly closed 
if and only if it is nilpotent.  


However, D. V. Osin noted that there exist finitely generated 
non-nilpotent unimodularly closed groups.
Indeed, let
$w_1(x)=1,\;w_2(x)=1,\dots$ be all unimodular equations over the
free group of rank two, let
$n_1,n_2,\dots$ be positive integers,
and let 
$\~x_i\gamma_{n_i}(F)$
be a solution to $w_i(x)=1$ in 
$F/\gamma_{n_i}(F)$ 
(existing by the Shmelkin theorem).
Then 
$w_i(\~x_i)\in\gamma_{n_i}(F)$,
 and the finitely generated unimodularly 
closed quotient group~$G=F/\nc{w_1(\~x_i),w_2(\~x_i),\dots}$ is a 
Golod--Shafarevich group if the integers~$n_i$ grow fast enough 
(see [E12]); in particular, $G$ is non-nilpotent.  
In the non-nilpotent finitely generated quotient group
$G/\bigcap\gamma_i(G)$, each unimodular equation is even uniquely 
solvable (the uniqueness follows from the Shmelkin theorem, since
this quotient group is residually nilpotent).

It is known that the following groups are 
also unimodularly (and even nonsingularly) closed:  
\- 
connected compact 
Lie groups [GR62]; 
\- 
Hall's universal locally finite group [Ha59] (its 
nonsingular closedness follows almost immediately from the results of 
[GR62]).

In Section 4, we prove inverse and implicit function theorems
implying, in particular, that
\disp{\sl\narrower\hfuzz 2cm
for any nilpotent group $G$, the
set of unimodular mappings $G^n\to G^n$
{\rm(which are defined naturally, see Section~4)}
forms a group under the composition.
}%
The structure of this group $\GA_n(G)$, which we call the
\emph{complete affine group}, is a subject for further study;
in the simplest case $G=\Z$, we
obtain the well-known group~$\GA_n(\Z)\iso\GL_n(\Z)\semitimes\Z^n$.

For non-nilpotent groups, unimodular mappings can be
very peculiar.
Anashin's results [An77] imply that,
for any finite non-abelian simple group $G$, all mapping
$G^n\to G^n$ are unimodular (for all $n$).
In Section 2, we construct
a surjective but non-injective unimodular mapping
from a group to itself.

Without unimodularity condition, very strange things can happen:
there exists an equation over an infinite finitely generated group
having
precisely one non-solution [KT05]
(in Section 2, we use a significantly stronger result
from [KOO13]).

\enditem
In relation to the Shmelkin theorem,
let us mention Kuzmin's result [Ku74]
(see also [Ku06]): {a finitely generated metabelian group $G$
is residually nilpotent if and only if
any unimodular equation
with one unknown
over $G$ has at most one solution in $G$.}
The same is true for finitely generated groups
with an abelian term of the lower central series [Ku78]).
Other recent results on nonsingular and unimodular
equations over groups can be found, e.g., in
[NT22],
[KM23],
[M24],
[KMR24],
[M26a],
[M26b],
[M26c].

In Section 5, we extend results from Section 3 to Lie algebras
and show that
\disp{\sl
a finitely generated solvable
Lie algebra \(over any field\) is nilpotent
if and only if
any
nonsingular equation over this algebra has a solution in this algebra.
}%
Here, a \emph{nonsingular equation} over a Lie algebra $A$
is the equality of the form $w(x)=0$, where~$w(x)$ is an element of the 
free product $A*\gp x$ not lying in the ideal generated by~$A$ (and 
$\gp x$ is the one-dimensional Lie algebra generated by $x$).

Actually, the Lie-algebra case is tricker:
for groups (in Section 3), we use that
any finitely generated
non-nilpotent solvable group has a finite non-nilpotent homomorphic
image;
while for Lie algebras, a similar assertion is false.
Indeed, the three-dimensional nilpotent Lie algebra
$$
L=\pres<a,b,c|[a,b]=c,\; [a,c]=[b,c]=0>
\qbox{(over a field $F$)}
$$
has an infinite-dimensional representation
$L\to \gl(F[t])$, where $a\mapsto(\hbox{multiplication by $t$})$,
$b\mapsto{{\rm d}\over {\rm d}t}$, $c\mapsto\id$.
This representation is irreducible if $\Char F=0$; therefore,
the corresponding semidirect sum of $L$ and the ideal $F[t]$
has all finite-dimensional images
nilpotent (because the ideal $F[t]$ is infinite-dimensional and
minimal).

\smallskip

{\noindent \bf Our notation}
is mainly standard. Note only that, if
$k\in \Z$ and $x$ and $y$ are elements of a group, then $x^y$, $x^{ky}$
and $x^{-y}$ denote $y^{-1}xy$, $y^{-1}x^ky$ and $y^{-1}x^{-1}y$,
respectively.
The commutator $[x,y]$ (in a group) is understood as $x^{-1}y^{-1}xy$.
The symbol~$\gp x_n$ denotes the cyclic group of order $n$ generated
by~$x$.
The free group of rank $n$ (with a basis $x_1,\dots,x_n$) is denoted
by $F_n$ (or $F(x_1,\dots,x_n)$). The symbol~$*$ denotes the free product,
and $\nc X$ is the normal closure of a set $X$
(in a group); $\gamma_i(G)$ stands for the $i$th term of the lower central
series of a group $G$. The $n$th direct power of a group~$G$ is denoted 
by~$G^n$.

\medskip

The authors thank Denis Osin for a valuable remark on finitely 
generated non-nilpotent unimodularly closed groups (see above).
The second author thanks
the Theoretical Physics and Mathematics Advancement Foundation ``BASIS".


\s 2.
Surjective but not injective unimodular mappings

A \emph{unimodular} mapping from a group $G$ to itself is a mapping
of the form $g\mapsto w(g)$, where the word $w\in G*\gp x_\infty$ is
unimodular.

The following result is a simplified
and weakened version of Theorem 2.3 from [KOO13].

\Th KOO {\rm[KOO13]}.
Any
countable group
$G_1$ without elements of order two
embeds into a simple 2-generated group $\KOO$, over which
there exists an equation $w(x)=1$ with precisely one
non-solution in $\KOO$ {\rm(i.e., the solutions to this equation are all
elements of $\KOO$, except for precisely one)}.

\Lemma 1.
There exist a torsion-free group $G_1$ and a
singular equation
$v(x)=1$ over $G_1$ such that, in any larger group
$\~G_1\supseteq G_1$, the set $M$ of solutions to this equation is such
that $bM\subsetneqq M$ for some $b\in G_1$.

\Proof
We put $G_1=\pres<a,b,c|a^b=a^2,[a^2,c]=1>$
(this is the free product
with amalgamated cyclic subgroups
of the Baumslag--Solitar group $\BS(1,2)$ and the free abelian group of
rank two) and $v(x)=[a^x,c]\in G_1*\gp x_\infty$.

If $\~x\in\~G_1$ is a solution
(i.e., $v(\~x)=[a^{\~x},c]=1$), then
$$
v(b\~x)=\[a^{b\~x},c\]=\[a^{2\~x},c\]=\[\(a^{\~x}\)^2,c\]=1.
$$
Thus, $bM\subseteq M$.
To show that the inclusion is strict, note that the identity element is not a
solution ($v(1)=[a,c]\ne1$),
and $b\cdot1$ is a solution: $v(b)=\[a^b,c\]=\[a^2,c\]=1$.
Thus, $b\in M\setminus bM$.
This completes the proof.

\Th 1.
There exist a 2-generated simple group~$G$ and a unimodular word
$u\in G*\gp x_\infty$ such that the mapping $G\to G$, $g\mapsto u(g)$ is
surjective but not injective.

\Proof
We put $G=\KOO$, where $G_1$ is the group from Lemma 1.
Let $w(x)=1$ be an equation over $G$ with precisely one non-solution.

\-
First, we may assume that the unique non-solution to the equation
$w(x)=1$ is the identity element
(this can be achieved by an obvious change of variables).

\-
Secondly, we may assume that $w(1)$ is an
arbitrary given
nonidentity element of $\KOO$.
Indeed, each element $g'\in G$ lies in the normal
closure of
$w(1)$ (since $G=\KOO$ is simple),
i.e., $g'=\prod w(1)^{\pm g_i}$. Then all
nonidentity elements of $G$ are solutions to
$w'(x)\equiv\prod w(x)^{\pm g_i}=1$, and $w'(1)=g'$.

\-
In particular,
there exists $w(x)\in G*\gp x_\infty$ such that
$w(g)=
\cases{
1&under $g\ne1$;
\cr
b^{-1}&under $g=1$;
}
$
(where~$b\in G_1$ is from Lemma 1).

\-
{%
Now, the unimodular word $u(x)=w\(v(x)\)x$ is as required.
Indeed, by Lemma 1 (applied to~$\~G_1=\KOO=G$)
\itemitem{--}
for $g\notin M$,
we have $u(g)=w(\underbrace{v(g)}_{\ne1})\cdot g=1\cdot g=g$;
\itemitem{--}
for $g\in M$,
we have $u(g)=w(\underbrace{v(g)}_{=1})\cdot g=b^{-1}\cdot g$.

Thus, the word $u(x)$ defines a mapping, which acts
identically outside the set $M$,
and maps the set $M$ onto a larger set
$b^{-1}M\supsetneqq M$ (by Lemma~1).
This completes the proof.
}


\s 3.
Residual finiteness of unimodular equations

\Th 2.
If all unimodular equation over a group $G$ have solutions in
a group $\~G\supseteq G$,
then any homomorphism from~$\~G$ to a finite group maps
$G$ to a nilpotent subgroup.

\Proof
The solvability of unimodular equations
is preserved under
homomorphic images; hence, we assume that the group $\~G$ is
finite, and we have to prove the nilpotency
of $G\subseteq\~G$.
First, Let us
show that
\dispno{\sl \narrower\hfuzz11cm
for any finite non-nilpotent
group~$G$, there exists a non-surjective
unimodular mapping $f\:G\to G$.
}(1)%
A hypothetical counterexample of minimal order
to this assertion
must have
a trivial centre, because the non-surjectivity of a unimodular mapping
$G/\ZZ(G)\to G/\ZZ(G)$ implies the non-surjectivity of
its (unimodular) lift $G\to G$,
and a group with a nilpotent central quotient is nilpotent.

In the nontrivial finite group $G$ with trivial centre,
we take
a minimal normal subgroup $B$,
an element $b\in B\setminus\1$,
and an element
$g\in G$ not commuting with $B$.
Then $B=\nc b=\nc{[g,b]}$, i.e.,
$b=\prod[g,b]^{\pm g_i}$, and the unimodular
equation~$x=\prod[g,x]^{\pm g_i}$ has two different solution: $x=b$ and
$x=1$, i.e.,
the unimodular mapping~$x\mapsto x\(\prod[g,x]^{\pm g_i}\)^{-1}$
is non-injective (and, therefore, non-surjective, because the group is
finite). This contradiction proves (1).

Proving from contradiction, suppose that $G$ is non-nilpotent.
Let us extend the
non-bijective
mapping $f$ from (1) to a unimodular mapping
$\~f\:\~G\to\~G$ (defined by the same
unimodular
word~$w$).
Note that
each power~$\~f^k\:=\~f\o\dots\o\~f$ of $\~f$
is unimodular too,
and, for some power~$F=\~f^k$, we have $F=F\o F$
(because each finite semigroup contains an idempotent).
The subgroup $G$ cannot belong to the image of~$F$,
because $F$ acts identically on this image,
while  $f$ and $F$ act non-injectively on $G$.
Therefore, the unimodular equation
$\underbrace{w(w(\dots(t)\dots))}_{k\hbox{ \small pairs of brackets}}=g$
(with coefficients from $G$) has no solutions in $\~G$ for
$g\in G\setminus F(\~G)$. This completes the proof.

\Corollary.
Suppose that a finitely generated group $G$ has one of the following
properties:
\item{\rm a)}
each nontrivial quotient of $G$ contains a nontrivial
finite or abelian normal subgroup
\(e.g., $G$ is solvable\)\;
\item{\rm b)}
$G$ is linear.
\enditem
Then $G$ is unimodularly closed
if and only if it is nilpotent.

\Proof
In a nilpotent group, all unimodular equations a solvable by the
Shmelkin theorem.
The other direction of the corollary follows
from Theorem~2
and a known fact:
\disp{\hfuzz 1cm
\sl\narrower\narrower\narrower\narrower\narrower
each finitely generated non-nilpotent group
with properties {\rm a)} or {\rm b)}
has a finite non-nilpotent image.
}%
The case a) was proven in [Ro70] (see also [Ro82]),
the case
b) was established in [Pl66] and [We68].

\s 4.
Inverse and implicit function theorems

Let $G$ be a group.
We consider
the \emph{formal polynomial mapping group}
$$
G(X)=\Gn\:=\bigl(G*F(x_1,\dots,x_n)\bigr)^n
\qbox{(where $X=\{x_1,\dots,x_n\}$ is a set of letters)}
$$
with the usual (component-wise) multiplication.
The components of a formal polynomial
mapping~$\Phi$ are denoted by~$\Phi_i$,
i.e., $\Phi=(\Phi_1,\dots,\Phi_n)$.
The group~$G(X)$ has an additional binary operation,
\emph{composition}:
$
(\Phi\o\Psi)_i=\Phi_i(\Psi)=\Phi_i(\Psi_1,\dots,\Psi_n),
$
under which $G(X)$ is a monoid with the unity~$e\:=(x_1,\dots,x_n)$.
The component-wise multiplication
and the composition are related by the one-sided distributivity:
$(\Phi\Psi)\o\Theta=(\Phi\o\Theta)(\Psi\o\Theta)$.
Another observation is that
\dispno{\sl
if $\Phi\in[G(X),G(X)]$, then
$\Phi\o(\Psi\Theta)\in(\Phi\o\Psi)\cdot[\nc{\Theta},G(X)]$.
}(2)%
Indeed,
the components $\Phi_i$ lie in the commutator subgroup of $\Gn$,
and the components~$\Theta_i$ are central
(modulo $[\nc{\Theta},G(X)]$), therefore,
$
\Phi_i(\Psi\Theta)=\Phi_i(\Psi).
$

\Lemma 2.
If
$\Psi\in G(X)$ and
$\Phi\in\gamma_k\bigl(G(X)\bigr)$,
where $k\ge2$,
then
$
(e\Psi)\o(e\Phi)\in e\Psi\Phi\cdot
\gamma_{k+1}\bigl(G(X)\bigr).
$

\Proof
Modulo $\gamma_{k+1}\bigl(G(X)\bigr)$, we have
$$
(e\Psi)\o(e\Phi)
=
\bigl(e\o(e\Phi)\bigr)\bigl(\Psi\o(e\Phi)\bigr)
=
e\Phi\bigl(\Psi\o(e\Phi)\bigr)
\=^2
e\Phi(\Psi\o e)
=
e\Phi\Psi
=
e\Psi\Phi,
$$
where the first equality is the distributivity.
This completes the proof.

\medskip

We call a formal polynomial mapping
$\Phi$
\emph{unimodular}
if the determinant of the matrix composed of the exponent sums of
$x_i$ in $\Phi_j$ is~$\pm1$.

Each formal polynomial mapping $\Phi$
defines a mapping $G^n\to G^n$ denoted by
$\Phi_G$.
A mapping~$G^n\to G^n$ of the form $\Phi_G$
for some $\Phi\in G(X)$ is called
\emph{polynomial};
if $\Phi$ is unimodular,
then the mapping~$\Phi_G$ is called
\emph{unimodular}.

By the Shmelkin theorem,
unimodular mappings $G^n\to G^n$ are bijective if
the group $G$ is nilpotent.

\proclaim Inverse function theorem.
If a group $G$ is
nilpotent, then the inverse to
any
unimodular mapping $\Gamma_G\:G^n\to G^n$
is unimodular too.

\Proof
We have to find a unimodular formal mapping
$\Psi\in G(X)$
such that the mapping $(\Psi\o\Gamma)_G$ is identical.

First, replacing $\Gamma$ with the composition $\Theta\o\Gamma$,
where $\Theta$ is a unimodular formal mapping, whose
matrix is inverse to the matrix of $\Gamma$,
we reduce the problem to the case, where the matrix of $\Gamma$
is identity.

Secondly,
replacing $\Gamma$ with the composition $\Xi\o\Gamma$,
where $\Xi=e\cdot\bigl(\Gamma_G(1,\dots,1)\bigr)^{-1}$,
we reduce the problem to the case, where
the matrix of $\Gamma$ is identity and
$\Gamma_G(1,\dots,1)=(1,\dots,1)$.

Thus, we assume that $\Gamma\in e\cdot\gamma_2\bigl(G(X)\bigr)$.
By Lemma 2, we have
$\Gamma'=\Delta'\o\Gamma\in e\cdot\gamma_3\bigl(G(X)\bigr)$
for some $\Delta'\in G(X)$
(namely, if $\Gamma=e\Phi$, then $\Delta'=e\Phi^{-1}$).
Applying Lemma~2 again and
proceeding in the same way,
we obtain
$
\Gamma^{(s)}=\Delta^{(s)}\o\Delta^{(s-1)}\o\dots\o\Delta'\o\Gamma
\in
e\cdot\gamma_{s+2}\bigl(G(X)\bigr),
$
where $s+2$ is larger than the nilpotency class of $G$. This means that
the mapping $\Gamma^{(s)}_G$ is identically and the 
mapping~$(\Delta^{(s)}\o\dots\o\Delta')_G$ is inverse to $\Gamma_G$.  This 
completes the proof.

\proclaim Implicit function theorem.
For any nilpotent group $G$ and
any system of equations
$$
S=\{
w_1(x_1,\dots,x_n,y_1,\dots,y_m)=1,
\dots,
w_n(x_1,\dots,x_n,y_1,\dots,y_m)=1
\}
$$
such that
the determinant of the
matrix composed of the exponent sums of unknowns $x_i$
in equations~$w_j=1$ is one,
there exist words
$
v_1(y_1,\dots,y_m),
\dots,
v_n(y_1,\dots,y_m)
\in G*F(y_1,\dots,y_m)
$
such
that the set of solutions to $S$
has the form
$$
\bigl\{
(\~x_1,\dots,\~x_n,\~y_1,\dots,\~y_m)\in G^{n+m}
\mid
\~x_1=v_1(\~y_1,\dots,\~y_m),
\dots,
\~x_n=v_n(\~y_1,\dots,\~y_m)
\bigr\}.
$$
\rm
Thus, $\~y_i$ are free variables,
i.e., they can take any values from~$G$.

\Proof
Suppose that
$
G_Y=\bigl(G*F(y_1,\dots,y_m)\bigr)/
\gamma_{s+1}\bigl(G*F(y_1,\dots,y_m)\bigr),
$
where
$s$ is the nilpotency class of $G$.
The words $w_i$ define a unimodular
mapping
$\Phi_{G_Y}\:G_Y^n\to G_Y^n$.
By the inverse function theorem, the inverse mapping
$(\Phi_{G_Y})^{-1}_\o$
is unimodular, i.e.,
(in particular)
for some words~$r_i(t_1,\dots,t_n)\in G_Y*F(t_1,\dots,t_n)$,
$$
(\Phi_{G_Y})^{-1}_\o(1,\dots,1)
=\bigl(
r_1(1,\dots,1),\dots,r_n(1,\dots,1)
\bigr).
$$
The elements $v_i=r_i(1,\dots,1)\in G_Y$ (considered as words in
$y_1,\dots,y_m$) are as required.
Indeed,
\-
the system of equations
$\{w_i(x_1,\dots,x_n,y_1,\dots,y_m)=1\}$
over $G_Y$ can be
rewritten as
$$
\Phi_{G_Y}(x_1,\dots,x_n)=(1,\dots,1);
$$
\-
$x_i=v_i$ is a solution to this system,
because,
substituting $x_i=v_i$
to the left-hand side, we obtain
$$
\Phi_{G_Y}(v_1,\dots,v_n)
=
\Phi_{G_Y}\bigl((\Phi_{G_Y})^{-1}_\o(1,\dots,1)\bigr)
=
(1,\dots,1);
$$
\-
each tuple $(g_1,\dots,g_m)\in G^m$ defines a retraction
$G_Y\to G$ ($y_i\mapsto g_i$);
\-
hence, each tuple
$\bigl(v_1(g_1,\dots,g_m),\dots,v_n(g_1,\dots,g_m),g_1,\dots,g_m\bigr)$
is a solution to $S$;
\-
conversely,
if we have a solution $x_i=g_i'$, $y_j=g_j$
to $S$ in~$G$, then the system of equations
$\{w_i(x_1,\dots,x_n,g_1,\dots,g_m)=1\}$ over $G$ has two solutions:
$\{x_i=g'\}$ and $\{x_i=v_i(g_1,\dots,g_m)\}$;
but a unimodular system over a nilpotent group
has a unique solution by the Shmelkin theorem;
therefore, $g_i'=v_i(g_1,\dots,g_m)$ as required.


\s 5.
Unimodular equations over Lie algebras

\proclaim Quillen theorem {\rm[Qu69] (see also
[Di78] or [Ba85])}.
If an associative algebra $U$ with unity over a field $F$ has
a filtration
$$
1\in U_0\subseteq U_1\subseteq\dots
\qbox{
\(where $U_i$ are subspace, $U_iU_j\subseteq U_{i+j}$
and $\bigcup U_i=U$\)}
$$
such that the associated graded algebra
$\gr(U)\:=\bigoplus U_i/U_{i-1}$
\(where $U_{-1}\:=\0$\)
is finitely generated and commutative,
then any endomorphism of any simple $U$-module is algebraic over $F$
\rm(i.e. it is annihilated by a nonzero polynomial from $F[x]$).

\Lemma 3.
If an associative algebra $U$ with unity over a field $F$
is equipped with a filtration
$$
1\in U_0\subseteq U_1\subseteq\dots
\qbox{
\(where $U_i$ are subspace, $U_iU_j\subseteq U_{i+j}$
and $\bigcup U_i=U$\)},
$$
and the polynomials algebra $U[t]$ is equipped with the filtration
$U[t]_i=\bigoplus\limits_{k+l=i}U_kt^l$,
then
the algebras~$\gr(U[t])$ and~$\gr(U)[t]$ are isomorphic
\(as abstract, not graded, algebras\).

\Proof
$
U[t]_i/U[t]_{i-1}=
\(\bigoplus\limits_{k+l=i}U_kt^l\)
\!\!\!\Bigm/\!\!\!
\(\bigoplus\limits_{p+q=i-1}U_pt^q\)
\iso
\bigoplus\limits_{k+l=i}(U_k/U_{k-1})t^l
$
(the last isomorphism is a general fact:
$\(\bigoplus C_s\)/\(\bigoplus D_s\)
\iso
\bigoplus(C_s/D_s),
$
if $D_s\subseteq C_s$).
Therefore,
$$
\gr(U[t])
=
\bigoplus U[t]_i/U[t]_{i-1}
=
\bigoplus\limits_i\bigoplus\limits_{k+l=i}(U_k/U_{k-1})t^l
=
\bigoplus\limits_{k,l}(U_k/U_{k-1})t^l.
$$
On the other hand,
$
\gr(U)[t]
=
\bigoplus\limits_{k}(U_k/U_{k-1})[t]
=
\bigoplus\limits_{k,l}(U_k/U_{k-1})t^l
$.
It is easy to show that the multiplication of elements of subspaces
$(U_k/U_{k-1})t^l$ in algebras $\gr(U[t])$ and $\gr(U)[t]$
is the same
(namely, the natural one:
$
(u_k+U_{k-1})t^l\cdot(u_j+U_{j-1})t^m
=(u_ku_j+U_{k+j-1})t^{l+m}
$). 
This completes the proof.

\Lemma 4.
Suppose that an associative algebra
$U$ over a field $F$
has an increasing filtration such that the associated
graded algebra is finitely generated and commutative,
and $\phi$ is
an endomorphism of finitely generated $U$-module $M$ such that
the endomorphism
$f(\phi)$ is surjective
for any polynomial $f\in F[x]$ with nonzero free term. Then
$\phi^n=0$ for some $n\in\N$.

\Proof
The surjectivity of the endomorphism $f(\phi)$
is inherited by quotient modules modulo $\phi$-invariant submodules.
Let us show that
\disp{\sl \hfuzz 1cm
among $\phi$-invariant
submodules $N\subset M$ not containing
$\phi^n(M)$ for any $n$,
there exists an inclusion-maximal submodule.
}%
By Zorn's lemma, it suffices to show that,
if the union
of a linearly ordered (by inclusion) family~$\cal N$
of
submodules contains $\phi^n(M)$ (for some $n$), then some
$N\in\cal N$
contains $\phi^n(M)$.
This follows immediately from the finite generatedness
of $M$: if $B$ is a finite set generating $M$,
then a module~$N\in\cal N$
contains the finite set $\phi^n(B)$ and, therefore,
$\phi^n(M)\subseteq N$.

Taking the quotient by a maximal $\phi$-invariant
submodule not containing
$\phi^n(M)$ for any $n$,
we obtain
(proving by contradiction) that the assertion of the lemma is false for
$M$, but true for any its proper quotient module by a $\phi$-invariant
submodule. Therefore, we assume that, for each nonzero $\phi$-invariant
submodule~$N\subseteq M$,
there exists $k\in\N$ such that $\phi^k(M)\subseteq N$. In particular,
the kernel $\ker\phi$ can be assumed to be zero (because otherwise
$\phi^k(M)\subseteq\ker\phi$ and $\phi^{k+1}=0$ as required).

Therefore,
$\phi$ is injective, and
$\phi(M)\iso M$.
Thus, we construct a bi-infinite chain of modules~%
$
\dots\subseteq M_{-1}\subseteq M_0\subseteq M_1\subseteq\dots
$,
where each module $M_i$ is isomorphic to $M$, and $\phi(M_i)=M_{i-1}$.
Moreover, $\phi$ is an automorphism of the
$U$-module $\=M=\bigcup M_i$.
($\phi$ extends to
the module $\=M$
naturally: $\phi(m_i)\:=\phi(m)_i$, where
$m\mapsto m_i$ is an isomorphism from $M$ to $M_i$.)
In addition,
\dispno{\sl
for any nonzero
$\phi$-invariant
submodule
$N\subseteq M_i$,
there exists $k\in\N$ such that $\phi^k(M_i)\subseteq N$,
}(3)%
because, $M$ has this property, as mentioned above.

The $U$-module $\=M$ can be considered as a $U[x,y]$-module,
where $x$ acts as $\phi$,
and $y$ acts as~$\phi^{-1}$.
Each nonzero submodule $N\subseteq\=M$
intersects some subspace $M_i$ nontrivially;
therefore, $M_{i-k}=\phi^k(M_i)\subseteqq^{3} N\cap M_i\subseteq N$
for some $k$, hence, all subspaces $M_s$ are contained in $N$,
because $M_s=y^{s-i+k}M_{i-k}$ for $s\ge i-k$. Thus, $N=\=M$, and
$U[x,y]$-module is simple.

By Lemma 3 (applied two times)
algebra $U[x,y]$ has a filtration such
that the associated graded algebra
is finitely generated and commutative
(because $U$ has such a
filtration). Therefore,
the Quillen theorem implies that the
automorphism~$\phi$ of $U[x,y]$-module $\=M$
is algebraic over~$F$.
The minimal polynomial of $\phi$
must be irreducible
(because the endomorphism algebra of a simple module is a skew field
by Schur's lemma).
This
contradicts the surjectivity of $f(\phi)$
for all polynomials~$f$ with nonzero free term
and completes the proof.

\Th 3.
In a finitely generated solvable Lie algebra $A$ over a field, each
nonsingular equation with one unknown
has a solution if and only if $A$ is nilpotent.

\Proof
In nilpotent algebras
(even not necessary Lie or associative)
each nonsingular equation
(and even each finite nonsingular system of equations)
has a solution [BO23].

To prove the other direction, note that the
solvability of nonsingular equations is inherited by quotient algebras;
therefore, the induction over the
 derived length
allows us to assume that the quotient algebra of $A$ by the
last nonzero term $P$ of the derived series is nilpotent. This means
that the ideal $P$ is abelian and has a finite codimension
(because any finitely generated nilpotent algebra is finite-dimensional).

Let $I$ be a nilpotent ideal in $A$ of minimal
codimension such that $I\supseteq P$.
We have to show that $I=A$.
Assuming the contrary, we take a nonzero element $z+I$ from centre of
$A/I$.
To obtain a contradiction,
it suffices to show that the ideal $I+Fz$ is nilpotent
(where $F$ is the main field).

A nonsingular (unimodular) equation
$$
x+c_1[z,x]+c_2[z,[z,x]]+\dots=b
\eqno{(*)}
$$
(with unknown $x$) has a solution for any $c_i\in F$
and any $b\in I$ (by the condition). This solution $\~x$ lies in $I$
(since $b\in I$ and $[z,\~x]\in I$ because of the centrality
of~$z+I\in A/I$).

The vector space $M=I/[I,I]$ is a module over the nilpotent
finite-dimensional algebra $A/I$.
In the module language,
the solvability of all equations $(*)$
means that, for any polynomial $f\in F[y]$ with free
term one, the
endomorphism
$f(\phi)$ of the module~$M$ is surjective, where $\phi\:m\mapsto[z,m]$
(this is an endomorphism, because $z+I$ is central in~$A/I$).

The universal enveloping algebra $U$ of the finite-dimensional algebra
$A/I$ has an increasing filtration such that the corresponding graded
algebra is the polynomial algebra in a finite number of variables (see,
e.g., [Ba85]). Therefore, Lemma 4 implies that $\phi^k=0$ for some
$k$. Hence, the ideal $(I+Fz)/[I,I]$ is nilpotent
of class at most $k$. This implies the
nilpotency of the ideal $I+Fz$ (and completes the proof) by virtue of
the following known analog of Hall's theorem.

\Th Chao--Stewart {\rm [Ch68], [St70]}.
If $I$ is a nilpotent ideal in a Lie algebra $L$, and the quotient algebra
$L/[I,I]$ is nilpotent, then $L$ is nilpotent.


\References

[An77]
V. S. Anashin,
Functionally complete groups,
Math. Notes, 22:1 (1977), 571-574.

[Ba85]
Yu. A. Bahturin,
Identities in Lie algebras.
M: Nauka, 1985.

[Di78]
G. Dixmier,
Enveloping algebras.
Amsterdam: North Holland, 1978.


[BO23]
Yu. Bahturin, A. Olshanskii,
Nilpotent algebras, implicit function theorem, and polynomial quasigroups,
Journal of Algebra 632 (2023), 154-193.
\arXiv:2208.09527

[Ch68]
Chong-Yun Chao,
Some characterizations of nilpotent Lie algebras,
Mathematische Zeitschrift, 103:1 (1968), 40-42.


[E12]
M. Ershov, 
Golod--Shafarevich groups: a survey,
International Journal of Algebra and Computation, 22:05 (2012), 1230001.
\arXiv:1206.0490

[GR62]
M. Gerstenhaber, O. S. Rothaus,
The solution of sets of equations in groups,
{Proc. Nat. Acad. Sci. USA}, {48:9} (1962), 1531-1533.

[FeR96]
R. Fenn, C. Rourke,
Klyachko's methods and the solution of
equations over torsion-free groups,
{L'Enseignment
Math\'ematique}, 42 (1996), 49-74.

[Ha59]
Ph. Hall,
Some Constructions for Locally Finite Groups,
Journal of the London Mathematical Society, s1-34:3 (1959),
305-319.


[K93]
A. A. Klyachko,
A funny property of sphere and equations over groups,
Communications in Algebra, 21:7 (1993), 2555-2575.

[KMR24]
A. A. Klyachko, M. A. Mikheenko, V. A. Roman'kov,
Equations over solvable groups,
Journal of Algebra, 638 (2024), 739-750.
\arXiv:2303.13240

[KOO13]
A. A. Klyachko, A. Yu. Olshanskii, D. V. Osin,
On topologizable and non-topologizable groups,
Topology and its Applications, 160:16 (2013), 2104-2120.
\arXiv:1210.7895

[KT05]
A. A. Klyachko, A. V. Trofimov,
The number of non-solutions of an equation in a group,
Journal of Group Theory, 8:6 (2005), 747-754.
\arXiv:math.GR/0411156

[Ku74]
Yu. V. Kuz'min,
Residual properties of metabelian groups [in Russian],
Algebra i Logika, 13:3 (1974), 300-310.

[Ku78]
Yu. V. Kuz'min,
Some approximation properties the variety $\eufm{AN}_c$,
Russian Math. Surveys, 33:4 (1978), 257-258.

[Ku06]
Yu. V. Kuz'min,
Homological group theory [in Russian].
Moscow: Faktorial, 2006.

[M24]
M. A. Mikheenko,
O $p$-nonsingular systems equations over solvable groups,
Sbornik: Mathematics,
215:6 (2024), 775-789.
\arXiv:2309.09096

[M26a]
M. A. Mikheenko,
Infinite systems of equations in abelian and nilpotent groups,
Journal of Group Theory, 29:1 (2026), 85-103.
\arXiv:2410.20729

[M26b]
M. A. Mikheenko,
Unimodular equations which do not preserve the derived length of a group,
Izv Ross. Akad. Nauk. Ser. Mat. (to appear).
\arXiv:2505.12783

[M26c]
M. A. Mikheenko,
On a generalization of Shmel'kin's theorem,
arXiv:2603.25110~.



[OAI26a]
OpenAI,
The Kervaire theorem for groups,
OpenAI Math Release preprint, 2026.

[OAI26b]
OpenAI,
Nonsingular systems of equations over arbitrary groups,
OpenAI Math Release preprint, 2026.


[Pl66]
V. P. Platonov,
Frattini subgroup of linear groups and residual finiteness [in Russian],
Dokl. Akad. Nauk SSSR, 171:4 (1966), 798-801.

[Qu69]
D. Quillen,
On the endomorphism ring of a simple module over an enveloping algebra,
Proc. Amer. Math. Soc, 21:1 (1969), 171-172.

[Ro70]
D. J. S. Robinson,
A theorem on finitely generated hyperabelian groups,
Inventiones mathematicae, 10:1 (1970), 38-43.

[Ro82]
D. J. S. Robinson,
A Course in the Theory of Groups.
New York: Springer, 1982.

[Sh67]
A. L. Shmel'kin,
Complete nilpotent groups [in Russian],
Algebra i Logika Seminar, 6:2 (1967), 111-114.

[St70]
I. N. Stewart,
Infinite-dimensional Lie algebras in the spirit
of infinite group theory,
Compositio Mathematica, 22:3 (1970), 313-331.

[We68]
B. A. F. Wehrfritz,
Frattini subgroups in finitely generated linear groups,
Journal of the London Mathematical Society, s1-43:1 (1968), 619-622.

\end